\documentclass[twoside, envcountsame, runningheads]{llncs}

\usepackage{epsfig}   
\usepackage{psfrag}
\usepackage{amsmath}
\usepackage{amsfonts}
\usepackage{amssymb}
\usepackage{rotating}
\usepackage{times}
\usepackage{bbm}
\usepackage{subfigure}
\usepackage{graphicx}

\usepackage{color}
\definecolor{nytt}{rgb}{0,0,1}
\definecolor{feil}{rgb}{1,0,0}
\definecolor{ut}{rgb}{1,.2,0}

\DeclareMathOperator{\Z}{Z}
\DeclareMathOperator{\I}{I}
\DeclareMathOperator{\length}{lg}

\newcommand\RR{\mathbb{R}}

\newcommand\CC{\mathbb{C}}

\renewcommand{\indexname}{\begin{tiny} Author \textit{Index} \end{tiny}}
 
\begin{document}

\newcommand{\kfield}{k}
\newcommand{\kbarfield}{\bar{k}}

\title{Monoid hypersurfaces}
\author{P{\aa}l~Hermunn~Johansen, Magnus~L{\o}berg, Ragni Piene}
\authorrunning{P. H. Johansen, M. L{\o}berg, R. Piene}
\institute{Centre of Mathematics for Applications and Department of Mathematics\\
 University of Oslo\\
 P. O. Box 1053 Blindern\\
 NO-0316 Oslo, Norway\\
\email{\{hermunn,mags,ragnip\}@math.uio.no}}

\maketitle

\begin{abstract}
A monoid hypersurface is an irreducible hypersurface of degree $d$ which has a
singular point of multiplicity $d-1$. Any monoid
hypersurface admits a rational parameterization, hence is of
potential interest in computer aided geometric design. We
study properties of monoids in general and of monoid surfaces
in particular. The main results include a description
of the possible real forms of the singularities on a monoid surface other
than the $(d-1)$-uple point. These results are applied to the
classification of singularities on quartic monoid surfaces,
complementing earlier work on the subject.
\end{abstract}

\section{Introduction}\label{sec_introduction}

A monoid hypersurface is an (affine or projective) irreducible algebraic 
hypersurface which has a singularity of multiplicity one 
less than the degree of the hypersurface. The presence of 
such a singular point forces the hypersurface to be 
rational: there is a rational parameterization given by (the 
inverse of) the linear projection of the hypersurface from 
the singular point.

The existence of an explicit rational parameterization makes 
such hypersurfaces potentially interesting objects in 
computer aided design. %
Moreover, since the ``space'' of 
monoids of a given degree is much smaller than the space of 
all hypersurfaces of that degree, one can hope to use 
monoids efficiently in (approximate or exact) 
implicitization problems.  These were the reasons for 
considering monoids in the paper 
\cite{DBLP:journals/cvgip/SederbergZKD99}. In 
\cite{MR2073068} monoid curves are used to approximate 
other curves that are close to a monoid curve, and in 
\cite{MR2116099} the same is done for monoid surfaces. 
In both articles the error of such approximations are analyzed 
-- for each approximation, a bound on the distance from the 
monoid to the original curve or surface can be computed.

In this article we shall study properties of monoid 
hypersurfaces and the classification of monoid surfaces with respect to their singularities. 
Section \ref{sec_basic_prop} explores properties of 
monoid hypersurfaces in arbritrary dimension and over an arbitrary base field. Section 
\ref{sec_monoid_surfaces} contains results on monoid 
surfaces, both over arbritrary fields and over $\RR$. 
The last section deals with the classification of monoid 
surfaces of degree four. Real and complex quartic monoid surfaces were first studied 
by Rohn \cite{MR1510280}, who gave a fairly complete description
of all possible cases. 
He also remarked \cite[p.~56]{MR1510280} that some of his results on 
quartic monoids hold for monoids of arbitrary degree; in particular, we believe he was aware
of many of the results in Section \ref{sec_monoid_surfaces}.
Takahashi, Watanabe, and Higuchi \cite{MR693454} classify \emph{complex} quartic 
monoid surfaces, but do not refer to Rohn. (They cite Jessop
\cite{Jessop1916}; Jessop, however, only treats quartic surfaces with double
points and refers to Rohn for the monoid case.)
Here we aim at giving a short description of the possible singularities that can 
occur on quartic monoids, with special emphasis on the real case.

\section{Basic properties}\label{sec_basic_prop}
Let $k$ be a field, let $\kbarfield$ denote its algebraic closure and 
$\mathbb{P}^n:=\mathbb{P}_{\kbarfield}^n$ the projective $n$-space
over $\bar{k}$. Furthermore we define the set of $\kfield$-rational
points $\mathbb{P}^n(\kfield)$ as the set of points that admit
representatives $(a_0:\dots:a_n)$ with each
$a_i\in\kfield$.

For any homogeneous polynomial $F\in\kbarfield[x_0,\dots,x_n]$ of
degree $d$ and point $p=(p_0:p_1:\dots:p_n)\in\mathbb{P}^n$ we can
define the multiplicity of $\Z(F)$ at $p$. We know that $p_r\neq0$ for
some $r$, so we can assume $p_0=1$ and write
\begin{equation*}
F=\sum_{i=0}^dx_0^{d-i}f_i(x_1-p_1x_0,x_2-p_2x_0,\dots,x_n-p_nx_0)
\end{equation*}
where $f_i$ is homogeneous of degree $i$. Then the multiplicity of $\Z(F)$
at $p$ is defined to be the smallest $i$ such that $f_i\neq0$.

Let $F\in\kbarfield[x_0,\dots,x_n]$ be of degree $d\ge3$. We say that
the hypersurface $X=\Z(F)\subset\mathbb{P}^n$ is a \emph{monoid}
hypersurface if $X$ is irreducible and has a singular point of
multiplicity $d-1$. 

In this article we shall only consider monoids $X=\Z(F)$ where the singular point
is $\kfield$-rational. Modulo a projective transformation of $\mathbb{P}^n$ 
over $\kfield$ we may -- and shall -- therefore assume that the singular point is the point
$O=(1:0:\cdots:0)$.

Hence, we shall from now on assume that $X=\Z(F)$, and
\begin{equation*}
F=x_0f_{d-1}+f_d,
\end{equation*}
where $f_i\in \kfield[x_1,\dots,x_n]\subset \kfield[x_0,\dots,x_n]$ is
homogeneous of degree $i$ and $f_{d-1}\neq0$. Since $F$ is
irreducible, $f_d$ is not identically 0, and $f_{d-1}$ and $f_d$ have
no common (non-constant) factors.

The \emph{natural rational parameterization} of the monoid $X=Z(F)$ is the map
\[\theta_F\colon \mathbb{P}^{n-1}\to \mathbb{P}^{n}\]
given by 
\[
\theta_F(a)=(f_{d}(a):-f_{d-1}(a)a_1:\ldots:-f_{d-1}(a)a_n),\]
for $a=(a_1:\dots:a_n)$ such that $f_{d-1}(a)\neq 0$ 
or $f_{d}(a)\neq 0$.

The set of lines through $O$ form a $\mathbb{P}^{n-1}$. For every
$a=(a_1:\dots:a_n)\in \mathbb{P}^{n-1}$, the line
\begin{equation}\label{eq_l_a}
L_{a}:=\{(s:ta_{1}:\ldots:ta_{n}) | (s:t)\in \mathbb{P}^1\}
\end{equation}
intersects $X=\Z(F)$ with multiplicity at least $d-1$ in $O$.
If $f_{d-1}(a)\neq 0$ or $f_{d}(a)\neq 0$, then the line $L_a$
also intersects $X$ in the point
\begin{equation*}
\theta_{F}(a)=(f_{d}(a):-f_{d-1}(a)a_1:\ldots:-f_{d-1}(a)a_n).
\end{equation*}
Hence the natural parameterization is the ``inverse" of the projection of $X$ from
the point $O$. Note that $\theta_F$ maps $\Z(f_{d-1})\setminus\Z(f_d)$ to $O$.  The
points where the parameterization map is not defined are called base
points, and these points are precisely the common zeros of $f_{d-1}$ and $f_{d}$.
Each such point $b$ corresponds to the line $L_b$ contained in the
monoid hypersurface. Additionally, every line of type $L_b$ contained
in the monoid hypersurface corresponds to a base point.

Note that $\Z(f_{d-1})\subset\mathbb{P}^{n-1}$ is the projective tangent cone to
$X$ at $O$, and that $\Z(f_d)$ is the intersection of $X$ with the
hyperplane ``at infinity'' $\Z(x_0)$.

Assume $P\in X$ is another singular point on the monoid $X$. Then the
line $L$ through $P$ and $O$ has intersection multiplicity at least
$d-1+2=d+1$ with $X$. Hence, according to Bezout's theorem, $L$ must
be contained in $X$, so that this is only possible if $\dim X\ge 2$.

By taking the partial derivatives of $F$ we can characterize the singular points of $X$ in terms of $f_d$ and $f_{d-1}$:

\begin{lemma}\label{lem_crit_sing}%
Let $\nabla=(\frac{\partial}{\partial
  x_1},\dots,\frac{\partial}{\partial x_n})$ be the gradient
operator.
\begin{itemize}
\item[$\mathrm{(i)}$]  A point $P=(p_0:p_1:\dots:p_n)\in\mathbb{P}^n$ is singular
  on $Z(F)$ if and only if $f_{d-1}(p_1,\dots,p_n)=0$ and $p_0\nabla
  f_{d-1}(p_1,\dots,p_n)+\nabla f_d(p_1,\dots,p_n)=0$. 
\item[$\mathrm{(ii)}$] All singular points of $\Z(F)$ are on lines $L_a$ where $a$ is a base point.
\item[$\mathrm{(iii)}$]  Both $\Z(f_{d-1})$ and $\Z(f_{d})$ are singular in a point $a\in
\mathbb{P}^{n-1}$ if and only if all points on $L_{a}$ are singular on
$X$.
\item[$\mathrm{(iv)}$]  If not all points on $L_{a}$ are singular, then at most one
point other than $O$ on $L_a$ is singular.
\end{itemize}
\end{lemma}
\begin{proof}
  (i) follows directly from taking the derivatives of
  $F=x_0f_{d-1}+f_d$, and (ii) follows from (i) and the fact that $F(P)=0$ for any singular point $P$. Furthermore, a point $(s:ta_{1}:\ldots:ta_{n})$ on $L_a$ is,
by (i),
  singular if and only if
\begin{equation*}
s\nabla f_{d-1}(ta)+\nabla f_d(ta)=t^{d-1}(s\nabla f_{d-1}(a)+t\nabla f_d(a))=0.
\end{equation*}
This holds for all $(s:t)\in\mathbb{P}^1$ if and only if $\nabla
f_{d-1}(a)=\nabla f_d(a)=0$. This proves (iii). If either $\nabla
f_{d-1}(a)$ or $\nabla f_{d-1}(a)$ are nonzero, the equation above has
at most one solution $(s_0:t_0)\in\mathbb{P}^1$ in addition to $t=0$,
and (iv) follows.
\end{proof}

Note that it is possible to construct monoids where $F\in\kfield[x_0,\dots,x_n]$,
but where no points of multiplicity $d-1$ are $\kfield$-rational. In
that case there must be (at least) two such points, and the line
connecting these will be of multiplicity $d-2$. Furthermore, the
natural parameterization will typically not induce a parameterization of
the $\kfield$-rational points from $\mathbb{P}^{n-1}(\kfield)$.

\section{Monoid surfaces}
\label{sec_monoid_surfaces}

In the case of a monoid surface, the parameterization has a finite
number of base points. 
From Lemma \ref{lem_crit_sing} (ii) we know that
all singularities of
the monoid other than $O$, are on lines $L_a$ corresponding to these
points. In what follows we will develop the theory for singularities
on monoid surfaces --- most of these results were probably known to Rohn
\cite[p.~56]{MR1510280}.

We start
by giving a precise definition of what we shall mean by a monoid surface.

\begin{definition}
\label{def_monoid_surface}
For an integer $d\ge3$ and a field $\kfield$ of characteristic $0$ the
polynomials $f_{d-1}\in\kfield[x_1,x_2,x_3]_{d-1}$ and
$f_d\in\kfield[x_1,x_2,x_3]_d$ define \emph{a normalized non-degenerate
  monoid surface} $\Z(F)\subset\mathbb{P}^3$, where
$F=x_0f_{d-1}+f_d\in\kfield[x_0,x_1,x_2,x_3]$ if the following hold:
\begin{itemize}
\item[$\mathrm{(i)}$]  $f_{d-1},f_d\neq0$
\item[$\mathrm{(ii)}$]  $\gcd(f_{d-1},f_d)=1$
\item[$\mathrm{(iii)}$]  The curves $\Z(f_{d-1})\subset\mathbb{P}^2$ and
  $\Z(f_d)\subset\mathbb{P}^2$ have no common singular
  point.
\end{itemize}
The curves $\Z(f_{d-1})\subset\mathbb{P}^2$ and
$\Z(f_d)\subset\mathbb{P}^2$ are called respectively the \emph{tangent
  cone} and the \emph{intersection with infinity}.
\end{definition}

Unless otherwise stated, 
a surface that satisfies
the conditions of Definition \ref{def_monoid_surface} shall be referred to simply
as a \emph{monoid surface}. 

Since we have finitely many base points $b$ and each line $L_b$
contains at most one singular point in addition to $O$, monoid surfaces
will have only finitely many singularities, so all singularities will
be isolated. (Note that Rohn included surfaces
with nonisolated singularities in his study \cite{MR1510280}.)
We will show that the singularities other than $O$ can be
classified by local intersection numbers.

\begin{definition}
  Let $f,g\in\kfield[x_1,x_2,x_3]$ be nonzero and homogeneous. Assume
  $p=(p_1:p_2:p_3)\in\Z(f,g)\subset\mathbb{P}^2$,
 and define the
  local intersection number
\begin{equation*}
\I_p(f,g)=\length\frac{\bar{\kfield}[x_1,x_2,x_3]_{m_p}}{(f,g)},
\end{equation*}
where $\bar\kfield$ is the algebraic closure of $\kfield$,
$m_p=(p_2x_1-p_1x_2,p_3x_1-p_1x_3,p_3x_2-p_2x_3)$ is the
homogeneous ideal of $p$, and $\length$ denotes the length
of the local ring
as a module over
itself.
\end{definition}

Note that $\I_p(f,g)\ge1$ if and only if $f(p)=g(p)=0$. When
$\I_p(f,g)=1$ we say that $f$ and $g$ intersect transversally at $p$.
The terminology is justified by the following lemma:

\begin{lemma}
  Let $f,g\in\kfield[x_1,x_2,x_3]$ be nonzero and homogeneous and
  $p\in\Z(f,g)$. Then the following are equivalent:
\begin{itemize}
\item[$\mathrm{(i)}$] $\I_p(f,g)>1$
\item[$\mathrm{(ii)}$]  $f$ is singular at $p$, $g$ is singular at $p$, or $\nabla f(p)$
and $\nabla g(p)$ are nonzero and parallel.
\item[$\mathrm{(iii)}$]  $s\nabla f(p)+t\nabla g(p)=0$ for some $(s,t)\neq(0,0)$
\end{itemize}
\end{lemma}
\begin{proof}
(ii) is equivalent to (iii) by a simple case study: $f$ is singular
at $p$ if and only if (iii) holds for $(s,t)=(1,0)$, 
$g$ is singular
at $p$ if and only if (iii) holds for $(s,t)=(0,1)$, 
and $\nabla
f(p)$ and $\nabla g(p)$ are nonzero and parallel if and only if (iii) holds for
some $s,t\neq 0$.

We can assume that $p=(0:0:1)$, so $\I_p(f,g)=\length S$ where
\begin{equation*}
S=\frac{\bar{\kfield}[x_1,x_2,x_3]_{(x_1,x_2)}}{(f,g)}.
\end{equation*}
Furthermore, let $d=\deg f$, $e=\deg g$ and write
\[f=\sum_{i=1}^df_ix_3^{d-i} \textrm{ and } g=\sum_{i=1}^eg_ix_3^{e-i}\] 
where $f_i,g_i$ are homogeneous of degree $i$.

If $f$ is singular at $p$, then $f_1=0$. Choose $\ell=ax_1+bx_2$ such
that $\ell$ is not a multiple of $g_1$. Then $\ell$ will be a nonzero
non-invertible element of $S$, so the length of $S$ is greater than
$1$.

We have $\nabla f(p)=(\nabla f_1(p),0)$ and $\nabla g(p)=(\nabla
g_1(p),0)$. If they are parallel, choose $\ell=ax_0+bx_1$ such that
$\ell$ is not a multiple of $f_1$ (or $g_1$), and argue as above.

Finally assume that $f$ and $g$ intersect transversally at $p$. We
may assume that $f_1=x_1$ and $g_1=x_2$.
Then $(f,g)=(x_1,x_2)$ as ideals in the local ring
$\bar{\kfield}[x_1,x_2,x_3]_{(x_1,x_2)}$. This means that $S$ is
isomorphic to the field $\bar{\kfield}(x_3)$. The length of any field
is $1$, so $\I_p(f,g)=\length S=1$.
\end{proof}

Now we can say which are the lines $L_b$, with $b\in\Z(f_{d-1},f_d)$, that
contain a singularity other than $O$:

\begin{lemma}
\label{lemma_criterion_sing_surface}
Let $f_{d-1}$ and $f_d$ be as in Definition \ref{def_monoid_surface}.
The line $L_b$ contains a singular point other than $O$ if and only if
$\Z(f_{d-1})$ is nonsingular at $b$ and the intersection
multiplicity $\I_b(f_{d-1},f_d)>1$.
\end{lemma}
\begin{proof}
Let $b=(b_1:b_2:b_3)$ and assume that $(b_0:b_1:b_2:b_3)$ is a
singular point of $\Z(F)$. Then, by Lemma \ref{lem_crit_sing},
$f_{d-1}(b_1,b_2,b_3)=f_d(b_1,b_2,b_3)=0$ and $b_0\nabla
f_{d-1}(b_1,b_2,b_3)+\nabla f_d(b_1,b_2,b_3)=0$, which implies
$\I_b(f_{d-1},f_d)>1$. Furthermore, if $f_{d-1}$ is singular at $b$,
then the gradient $\nabla f_{d-1}(b_1,b_2,b_3)=0$, so $f_d$, too, is
singular at $b$, contrary to our assumptions.
  
Now assume that $\Z(f_{d-1})$ is nonsingular at $b=(b_1:b_2:b_3)$ and
the intersection multiplicity $\I_b(f_{d-1},f_d)>1$. The second
assumption implies $f_{d-1}(b_1,b_2,b_3)=f_d(b_1,b_2,b_3)=0$ and
$s\nabla f_{d-1}(b_1,b_2,b_3)=t\nabla f_d(b_1,b_2,b_3)$ for some
$(s,t)\neq(0,0)$. Since $\Z(f_{d-1})$ is nonsingular at $b$, we know
that $\nabla f_{d-1}(b_1,b_2,b_3)\neq0$, so $t\neq0$. Now
$(-s/t:b_1:b_2:b_3)\neq(1:0:0:0)$ is a singular point of $\Z(F)$ on
the line $L_b$.
\end{proof}

Recall that an $A_n$ singularity is a singularity with normal form 
$x_1^2+x_2^2+x_3^{n+1}$, see \cite[p.~184]{MR777682}.

\begin{proposition}\label{prop_type_sing}
  Let $f_{d-1}$ and $f_d$ be as in Definition
  \ref{def_monoid_surface}, and assume $P=(p_0:p_1:p_2:p_3)\neq(1:0:0:0)$ is a
  singular point of $\Z(F)$ with $I_{(p_1:p_2:p_3)}(f_{d-1},f_d)=m$.
  Then $P$ is an $A_{m-1}$ singularity.
  \end{proposition}

\begin{proof}
  We may assume that $P=(0:0:0:1)$ and write the local equation
\begin{equation}
g:=F(x_0,x_1,x_2,1)=x_0f_{d-1}(x_1,x_2,1)+f_d(x_1,x_2,1)=\sum_{i=2}^dg_i
\end{equation}
with $g_i\in\bar{\kfield}[x_0,x_1,x_2]$ homogeneous of degree $i$.  Since
$\Z(f_{d-1})$ is nonsingular at $0:=(0:0:1)$, we can assume that the
linear term of $f_{d-1}(x_1,x_2,1)$ is equal to $x_1$. The quadratic
term $g_2$ of $g$ is then $g_2=x_0x_1+ax_1^2+bx_1x_2+cx_2^2$ for some
$a,b,c\in\kfield$. The Hessian matrix of $g$ evaluated at $P$ is
\begin{equation*}
H(g)(0,0,0)=H(g_2)(0,0,0)
=\left(\begin{matrix}0&1&0\\1&2a&b\\0&b&2c\end{matrix}\right)
\end{equation*}
which has corank $0$ when $c\neq0$ and corank $1$ when $c=0$. By
\cite[p.~188]{MR777682}, $P$ is an $A_1$ singularity when $c\neq0$ and an
$A_n$ singularity for some $n$ when $c=0$.

The index $n$ of the singularity is equal to the Milnor number
\begin{equation*}
\mu=\dim_{\bar\kfield}\frac{\bar\kfield[x_0,x_1,x_2]_{(x_0,x_1,x_2)}}{J_g}=\dim_{\bar\kfield}\frac{\bar\kfield[x_0,x_1,x_2]_{(x_0,x_1,x_2)}}{\left(\frac{\partial g}{\partial x_0},\frac{\partial g}{\partial x_1},\frac{\partial g}{\partial x_2}\right)}.
\end{equation*}
We need to show that $\mu=\I_0(f_{d-1},f_d)-1$. From the definition of
the intersection multiplicity, it is not hard to see that
\begin{equation*}
\I_0(f_{d-1},f_d)=\dim_{\bar\kfield}\frac{\bar\kfield[x_1,x_2]_{(x_1,x_2)}}{\left(f_{d-1}(x_1,x_2,1),f_d(x_1,x_2,1)\right)}.
\end{equation*}
The singularity at $p$ is isolated, so the Milnor number is finite.
Furthermore, since $\gcd(f_{d-1},f_d)=1$, the intersection multiplicity
is finite. Therefore both dimensions can be calculated in the
completion rings. For the rest of the proof we view $f_{d-1}$
and $f_d$ as elements of the power series rings
$\bar\kfield[[x_1,x_2]]\subset\bar\kfield[[x_0,x_1,x_2]]$, and all
calculations are done in these rings.

Since $\Z(f_{d-1})$ is smooth at $O$, we can write
\[
f_{d-1}(x_1,x_2,1)=\left(x_1-\phi(x_2)\right)u(x_1,x_2)\] 
for some
power series $\phi(x_2)$ and invertible power series $u(x_1,x_2)$. To
simplify notation we write $u=u(x_1,x_2)\in\bar\kfield[[x_1,x_2]]$.

The Jacobian ideal $J_g$ is generated by the three partial derivatives:
\begin{eqnarray*}
\frac{\partial g}{\partial x_0}&=&\left(x_1-\phi(x_2)\right)u \\
\frac{\partial g}{\partial x_1}&=&x_0\left(u+(x_1-\phi(x_2))\frac{\partial u}
{\partial x_1}\right)+\frac{\partial f_d}{\partial x_1}(x_1,x_2)\\
\frac{\partial g}{\partial x_2}&=&x_0\left(-\phi^\prime(x_2)u+(x_1-\phi(x_2))
\frac{\partial u}{\partial x_2}\right)
+\frac{\partial f_d}{\partial x_2}(x_1,x_2)
\end{eqnarray*}
By using the fact that $x_1-\phi(x_2)\in\left(\frac{\partial g}{\partial x_0}\right)$ we can write $J_g$ without the symbols $\frac{\partial u}{\partial x_1}$ and $\frac{\partial u}{\partial x_2}$:
\begin{equation*}
\textstyle
J_g=\left(x_1-\phi(x_2),x_0u+\frac{\partial f_d}{\partial x_1}(x_1,x_2),-x_0u\phi^\prime(x_2)+\frac{\partial f_d}{\partial x_2}(x_1,x_2)\right)
\end{equation*}

To make the following calculations clear, define the polynomials $h_i$
by writing $f_d(x_1,x_2,1)=\sum_{i=0}^dx_1^ih_i(x_2)$. Now
\begin{equation*}\textstyle
J_g=\left(x_1-\phi(x_2),x_0u+\sum_{i=1}^dix_1^{i-1}h_i(x_2),-x_0u\phi^\prime(x_2)+
\sum_{i=0}^dx_1^ih_i^\prime(x_2)\right),
\end{equation*}
so
\begin{equation*}
\frac{\bar\kfield[[x_0,x_1,x_2]]}{J_g}=\frac{\bar\kfield[[x_2]]}{(A(x_2))}
\end{equation*}
where
\begin{equation*}
\textstyle
A(x_2)=\phi^\prime(x_2)\left(\sum_{i=1}^di\phi(x_2)^{i-1}h_i(x_2)\right)+\left(\sum_{i=0}^d\phi(x_2)^ih_i^\prime(x_2)\right).
\end{equation*}
For the intersection multiplicity we have
\begin{equation*}
\frac{\bar\kfield[[x_1,x_2]]}{\Big(f_{d-1}(x_1,x_2,1),f_d(x_1,x_2,1)\Big)}
=\frac{\bar\kfield[[x_1,x_2]]}{\left(x_1-\phi(x_2),\sum_{i=0}^dx_1^ih_i(x_2)\right)}=\frac{\bar\kfield[[x_2]]}{\Big(B(x_2)\Big)}
\end{equation*}
where $B(x_2)=\sum_{i=0}^d\phi(x_2)^ih_i(x_2)$. Observing that
$B^\prime(x_2)=A(x_2)$ gives the result $\mu=\I_0(f_{d-1},f_d)-1$.
\end{proof}

\begin{corollary}\label{coro_max_sing}
  A monoid surface of degree $d$ can have at most $\frac{1}{2}d(d-1)$
  singularities in addition to $O$. If this number of singularities is
  obtained, then all of them will be of type $A_1$.
  \end{corollary}
\begin{proof}
  The sum of all local intersection numbers $\I_a(f_{d-1},f_d)$ is
  given by B\'ezout's theorem:
\begin{equation*}
\sum_{a\in\Z(f_{d-1},f_d)}\!\!\!\!\!\!\I_a(f_{d-1},f_d)=d(d-1).
\end{equation*}
The line $L_a$ will contain a singularity other than $O$ only if
$\I_a(f_{d-1},f_d)\ge2$, giving a maximum of $\frac{1}{2}d(d-1)$
singularities in addition to $O$. Also, if this number is obtained, all
local intersection numbers must be exactly $2$, so all singularities
other than $O$ will be of type $A_1$.
\end{proof}

Both Proposition \ref{prop_type_sing} and Corollary \ref{coro_max_sing}
were known to
Rohn, who stated these results only in the case $d=4$, but said they  could be
generalized to arbitrary $d$ \cite[p.~60]{MR1510280}.
\medskip

For the rest of the section we will assume $\kfield=\RR$. It
turns out that we can find a \emph{real} normal form for the singularities
other than $O$. The complex singularities of type $A_n$ come in several
real types, with normal forms $x_1^2\pm x_2^2\pm x_3^{n+1}$. Varying
the $\pm$ gives two types for $n=1$ and $n$ even, and three types for
$n\ge3$ odd. The real type with normal form $x_1^2-x_2^2+x_3^{n+1}$ is
called an $A_n^-$ singularity, or \emph{of type $A^{-}$}, and is what we find on real monoids:

\begin{proposition}\label{realsing} 
On a real monoid, all singularities other than $O$ are of type 
  $A^-$.
  \end{proposition}
\begin{proof}
  Assume $p=(0:0:1)$ is a singular point on $\Z(F)$ and set $g=F(x_{0},x_{1},x_{2},1)$
 as in the proof of Proposition \ref{prop_type_sing}.
  
  First note that $u^{-1}g=x_0(x_1-\phi(x_2))+f_d(x_1,x_2)u^{-1}$ is
  an equation for the singularity. We will now prove that $u^{-1}g$ is
  right equivalent to $\pm(x_0^2-x_1^2+x_2^n)$, for some $n$, by constructing right
  equivalent functions $u^{-1}g=:g_{(0)}\sim g_{(1)}\sim g_{(2)}\sim
  g_{(3)}\sim\pm(x_0^2-x_1^2+x_2^n)$. Let
\begin{eqnarray*}
g_{(1)}(x_0,x_1,x_2)&=&g_{(0)}(x_0,x_1+\phi(x_2),x_2)\\
&=&x_0x_1+f_d(x_1+\phi(x_2),x_2)u^{-1}(x_1+\phi(x_2),x_2)\\
&=&x_0x_1+\psi(x_1,x_2)
\end{eqnarray*}
where $\psi(x_1,x_2)\in\RR[[x_1,x_2]]$. Write
$\psi(x_1,x_2)=x_1\psi_1(x_1,x_2)+\psi_2(x_2)$ and define
\begin{equation*}
g_{(2)}(x_0,x_1,x_2)=g_{(1)}(x_0-\psi_1(x_1,x_2),x_1,x_2)=x_0x_1+\psi_2(x_2).
\end{equation*}
The power series $\psi_2(x_2)$ can be written on the form
\begin{equation*}
\psi_2(x_2)=sx_2^n(a_0+a_1x_2+a_2x_2^2+\dots)
\end{equation*}
where $s=\pm1$ and $a_0>0$. We see that $g_{(2)}$ is right equivalent
to $g_{(3)}=x_0x_1+sx_2^n$ since
\begin{equation*}
g_{(2)}(x_0,x_1,x_2)=
g_{(3)}\left(x_0,x_1,x_2\sqrt[n]{a_0+a_1x_2+a_2x_2^2+\dots}\right).
\end{equation*}
Finally we see that
\begin{equation*}
g_{(4)}(x_0,x_1,x_2):=g_{(3)}(sx_0-sx_1,x_0+x_1,x_2)=s(x_0^2-x_1^2+x_2^n)
\end{equation*}
proves that $u^{-1}g$ is right equivalent to $s(x_0^2-x_1^2+x_2^n)$
which is an equation for an $A_{n-1}$ singularity with normal form
$x_0^2-x_1^2+x_2^n$.
\end{proof}


Note that for $d=3$, the singularity at $O$ can be an $A_1^{+}$ singularity.
This happens for example when
$f_2=x_0^2+x_1^2+x_2^2$.
\medskip

\begin{figure}[!htbp]
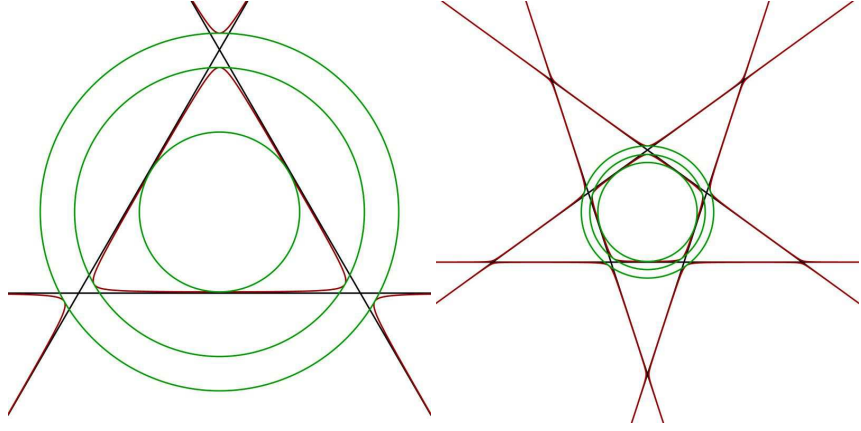

\begin{center}
  \includegraphics[width=.46\textwidth]{fig3_800_8.epsf}
  \includegraphics[width=.46\textwidth]{fig5_800_8.epsf}
\caption{The curves $f_m$ for $m=3,5$ and corresponding circles}\label{fig_real_sing}
\end{center}
\end{figure}

For a \emph{real} monoid, Corollary \ref{coro_max_sing} implies that we can
have at most $\frac{1}{2}d(d-1)$ \emph{real} singularities in addition to $O$.
We can show that the bound is sharp by a simple construction:
\smallskip

\noindent\textbf{Example.} 
To construct a monoid with the maximal number of real singularities,
it is sufficient to construct two \emph{affine} real curves
in the $xy$-plane defined by equations
$f_{d-1}$ and $f_d$ of degrees $d-1$ and $d$ such that the curves intersect in
$d(d-1)/2$ points with multiplicity $2$. Let $m\in\{d-1,d\}$ be odd
and set
\begin{equation*}
f_m=\varepsilon-\prod_{i=1}^m\left(x\sin\left(\frac{2i\pi}{m}\right)+y\cos\left(\frac{2i\pi}{m}\right)+1\right).
\end{equation*}
For $\varepsilon>0$ sufficiently small there exist at least
$\frac{m+1}{2}$ radii $r>0$, one for each root of the univariate polynomial $f_m\vert_{x=0}$, such that the circle $x^2+y^2-r^2$
intersects $f_m$ in $m$ points with multiplicity $2$. Let $f_{2d-1-m}$
be a product of such circles. Now the homogenizations of $f_{d-1}$ and
$f_d$ define a monoid surface with $1+\frac{1}{2}d(d-1)$
singularities. See Figure \ref{fig_real_sing}.\medskip
 
\begin{figure}[!htbp]
\begin{center}
  \includegraphics[width=.46\textwidth]{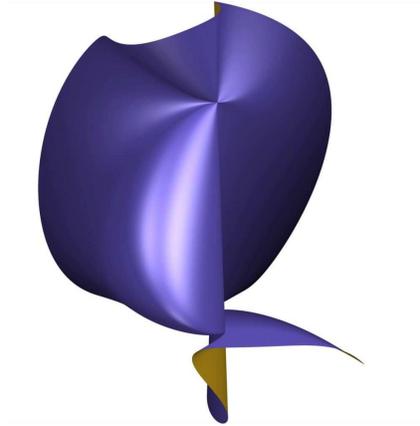}
\caption{The surface defined by $F=x_0(x_1x_2^{d-2}+x_3^{d-1})+x_1^d$ for $d=4$.}
\label{fig_ex_2}
\end{center}
\end{figure}

Proposition \ref{prop_type_sing} and Bezout's theorem imply that the maximal 
Milnor number of a singularity other than $O$ is $d(d-1)-1$. 
The following
example shows that this bound can be achieved on a real monoid:
\smallskip

\noindent\textbf{Example.} The surface $X\subset\mathbb{P}^3$ defined by 
$F=x_0(x_1x_2^{d-2}+x_3^{d-1})+x_1^d$ has exactly two singular points. The point 
$(1:0:0:0)$ is a singularity of multiplicity $d-1$
with Milnor number $\mu=(d^2-3d+1)(d-2)$,
while the point $(0:0:1:0)$ is an $A_{d(d-1)-1}$ singularity. 
A picture of this surface for
$d=4$ is shown in Figure \ref{fig_ex_2}.

\section{Quartic monoid surfaces}

Every cubic surface with isolated singularities is a monoid. Both smooth and
singular cubic surfaces have been studied extensively, most notably in
\cite{Schlafli1863}, where real cubic surfaces and their singularities
were classified, and more recently in \cite{MR0008171}, \cite{MR533323},
and \cite{MR888126}. The site \cite{LvSCSH} contains
additional pictures and references.

In this section we shall consider the case $d=4$.
The classification of real and complex \emph{quartic} monoid surfaces was started by Rohn
\cite{MR1510280}. (In addition to considering the singularities, Rohn studied
the existence of lines \emph{not} passing through the triple point, and that of
other special curves on the monoid.) In
\cite{MR693454}, Takahashi, Watanabe, and Higuchi
 described
the singularities of such \emph{complex} surfaces.
The monoid singularity  of a quartic monoid is minimally elliptic
\cite{MR646044}, and minimally elliptic singularities have the same complex
topological type if and only if their dual graphs are isomorphic
\cite{MR0568898}. In \cite{MR0568898} all possible dual
graphs for minimally elliptic singularities are listed, along with example
equations.

Using Arnold's notation for the singularities, we use and extend the approach
of Takahashi, Watanabe, and Higuchi in \cite{MR693454}.

Consider a quartic monoid surface, $X=\Z(F)$, with $F=x_{0}f_{3}+f_{4}$.
The tangent cone, $\Z(f_3)$, can be of one of nine (complex) types,
each needing a separate analysis.

For each type we fix $f_3$, but any
other tangent cone of the same type will be projectively equivalent
(over the complex numbers) to this fixed $f_3$. The nine different
types are:
\begin{enumerate}
\item Nodal irreducible curve, $f_3=x_1x_2x_3+x_2^3+x_3^3$.
\item Cuspidal curve, $f_3=x_1^3-x_2^2x_3$.
\item 
Conic and a chord, $f_3=x_3(x_1x_2+x_3^2)$
\item Conic and a tangent line, $f_3=x_3(x_1x_3+x_2^2)$.
\item 
Three general lines, $f_3=x_1x_2x_3$.
\item 
Three lines meeting in a point, $f_3=x_2^3-x_2x_3^2$
\item A double line and another line, $f_3=x_2x_3^2$
\item A triple line $f_3=x_3^3$
\item A smooth curve, $f_3=x_1^3+x_2^3+x_3^3+3ax_0x_1x_3$ where $a^3\neq-1$
\end{enumerate}

To each quartic monoid we can associate, in addition to the type,
several integer invariants, all given as intersection numbers. From
\cite{MR693454} we know that, for the types 1--3,
5, and 9, these invariants will determine the singularity type
of $O$ up to right equivalence. In
the other cases the singularity series, as defined by Arnol'd in
\cite{ArnoldLondon1} and \cite{ArnoldLondon2}, is determined by the
type of $f_3$. We shall use, without proof, the results on the singularity type of $O$
due to \cite{MR693454}; however, we shall use the notations
of \cite{ArnoldLondon1} and \cite{ArnoldLondon2}.
\bigskip

We complete the classification begun in \cite{MR693454} by supplying a complete
list of the possible singularities occurring on a quartic monoid. In addition, we extend
the results to the case of \emph{real} monoids. Our results are
summarized in the following theorem. 



\begin{theorem} On a quartic monoid surface, singularities other than the monoid point 
can occur as given in Table \ref{table_class_sing_quartic_monoids}. Moreover, all possibilities
are realizable on \emph{real} quartic monoids with a real monoid point, and with the other 
singularities being real and of type $A^{-}$.
\end{theorem}

\begin{table}[!htbp]
\begin{tabular}{c|c|l|l}
Case & Triple point & Invariants and constraints & Other singularities\\\hline
$1$ & $T_{3,3,4}$ && $A_{m_i-1}$, $\sum m_i=12$\\\cline{2-4}
& $T_{3,3,3+m}$ & $m=2,\dots,12$ & $A_{m_i-1}$, $\sum m_i=12-m$\\\hline
2 & $Q_{10}$ $$ && $A_{m_i-1}$, $\sum m_i=12$\\\cline{2-4}
& $T_{9+m}$ & $m=2,3$ & $A_{m_i-1}$, $\sum m_i=12-m$\\\hline
3 & $T_{3,4+r_0,4+r_1}$ & $r_0=\max(j_0,k_0)$, $r_1=\max(j_1,k_1)$, & $A_{m_i-1}$, $\sum m_i=4-k_0-k_1$,\\
&& $j_0>0\leftrightarrow k_0>0$, $\min(j_0,k_0)\le1$, &  $A_{m_i^{\prime}-1}$, $\sum m_i^{\prime}=8-j_0-j_1$\\
&&  $j_1>0\leftrightarrow k_1>0$, $\min(j_1,k_1)\le1$\\\hline
4 & $S$ series & $j_0\le8$, $k_0\le4$, $\min(j_0,k_0)\le2$, & $A_{m_i-1}$, $\sum m_i=4-k_0$,\\
&& $j_0>0\leftrightarrow k_0>0$, $j_1>0\leftrightarrow k_0>1$ & $A_{m_i^{\prime}-1}$, $\sum m_i^{\prime}=8-j_0$\\\hline
5&$T_{4+j_k,4+j_l,4+j_m}$ & $m_1+l_1\le4$, $k_2+m_2\le4$, & $A_{m_i-1}$, $\sum m_i=4-m_1-l_1$,\\
&& $k_3+l_3\le4$, $k_2>0\leftrightarrow k_3>0$, & $A_{m_i^{\prime}-1}$, $\sum m_i^{\prime}=4-k_2-m_2$,\\
&&  $l_1>0\leftrightarrow l_3>0$, $m_1>0\leftrightarrow m_2>0$, & $A_{m_i^{\prime\prime}-1}$, $\sum m_i^{\prime\prime}=4-k_3-l_3$\\
&& $\min(k_2,k_3)\le1$, $\min(l_1,l_3)\le1$,\\
&& $\min(m_1,m_2)\le1$, $j_k=\max(k_2,k_3)$,\\
&& $j_l=\max(l_1,l_3)$, $j_m=\max(m_1,m_2)$\\\hline
6 & $U$ series & $j_1>0\leftrightarrow j_2>0\leftrightarrow j_3>0$, &$A_{m_i-1}$, $\sum m_i=4-j_1$,\\
&& at most one of $j_1,j_2,j_3>1$, & $A_{m_i^{\prime}-1}$, $\sum m_i^{\prime}=4-j_2$,\\
&& $j_1,j_2,j_3\le4$ & $A_{m_i^{\prime\prime}-1}$, $\sum m_i^{\prime\prime}=4-j_3$\\\hline
7 & $V$ series & $j_0>0\leftrightarrow k_0>0$, $\min(j_0,k_0)\le1$, & $A_{m_i-1}$, $\sum m_i=4-j_0$,\\
&& $j_0\le4$, $k_0\le4$\\\hline
8 & $V^\prime$ series && None\\\hline
9 & $P_8=T_{3,3,3}$ && $A_{m_i-1}$, $\sum m_i=12$\\\hline
\end{tabular}
\smallskip

\caption{Possible configurations of singularities for each case}\label{table_class_sing_quartic_monoids}
\end{table}

\begin{proof}
The invariants listed in the ``Invariants and constraints'' column are
all nonnegative integers, and any set of integer values satisfying the
equations represents one possible set of invariants, as described
above. Then, for each set of invariants, (positive) intersection
multiplicities, denoted $m_i$, $m_i^{\prime}$ and $m_i^{\prime\prime}$,
will determine the singularities other than $O$. The column ``Other
singularities'' give these and the equations they must satisfy. Here
we use the notation $A_0$ for a line $L_a$ on $\Z(F)$ where $O$ is the
only singular point.

The analyses of the nine cases share many similarities, and we have
chosen not to go into great detail when one aspect of a case differs
little from the previous one.  We end the section with a discussion on
the possible real forms of the tangent cone and how this affects the
classification of the \emph{real} quartic monoids.

In all cases, we shall write
\begin{eqnarray*}
f_4&=&a_{1}x_1^4+a_{2}x_1^3x_2+a_{3}x_1^3x_3+a_{4}x_1^2x_2^2+a_{5}x_1^2x_2x_3\\
&+&a_{6}x_1^2x_3^2+a_{7}x_1x_2^3+a_{8}x_1x_2^2x_3+a_{9}x_1x_2x_3^2+a_{10}x_1x_3^3\\
&+&a_{11}x_2^4+a_{12}x_2^3x_3+a_{13}x_2^2x_3^2+a_{14}x_2x_3^3+a_{15}x_3^4
\end{eqnarray*}
and we shall investigate how the coefficients $a_1,\dots,a_{15}$ are related
to the geometry of the monoid.
\smallskip

{\bf Case 1.} The tangent cone is a nodal irreducible curve, and
we can assume 
\[
f_3(x_1,x_2,x_3)=x_1x_2x_3+x_2^3+x_3^3.\]
The nodal
curve is singular at $(1:0:0)$. If $f_4(1,0,0)\neq0$, then $O$ is a
$T_{3,3,4}$ singularity \cite{MR693454}. We recall that $(1:0:0)$ cannot be a singular
point on $\Z(f_4)$ as this would imply a singular line on the monoid,
so we assume that either $(1:0:0)\not\in\Z(f_4)$ or $(1:0:0)$ is a
smooth point on $\Z(f_4)$.  Let $m$ denote the intersection number
$I_{(1:0:0)}(f_3,f_4)$. Since $\Z(f_3)$ is singular at $(1:0:0)$ we
have $m\neq1$. 
From \cite{MR693454} we know that $O$ is a $T_{3,3,3+m}$ singularity for
$m=2,\dots,12$. Note that some of these complex singularities have two real forms,
as illustrated in Figure \ref{fig_T335}.

\begin{figure}[!htbp]
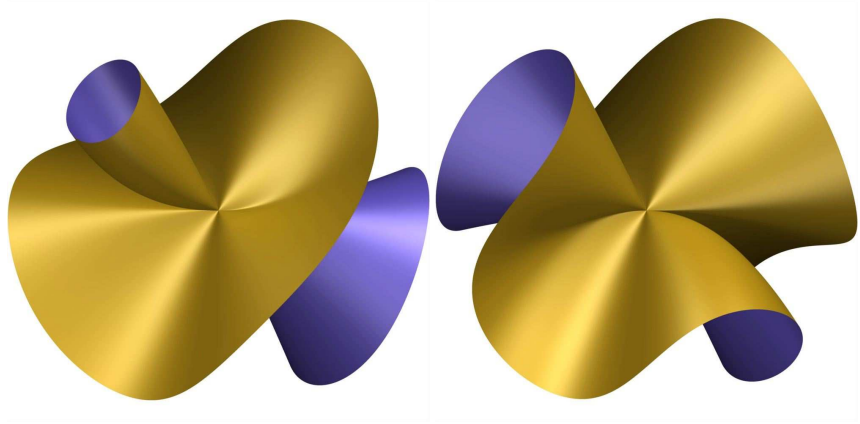
%
\includegraphics[width=.46\textwidth]{T-335_plus_scale0_2.epsf}
\includegraphics[width=.46\textwidth]{T-335_minus_scale0_2.epsf}%
\caption{The monoids $\Z(x^3+y^3+5xyz-z^3(x+y))$ and
$\Z(x^3+y^3+5xyz-z^3(x-y))$ both have a $T_{3,3,5}$ singularity, but the
singularities are not right equivalent over $\RR$. (The pictures are generated
by the program \cite{surf}.)}
\label{fig_T335}%
\end{figure}

B\'ezout's theorem and Proposition \ref{prop_type_sing} limit the
possible configurations of singularities on the monoid for each $m$.
Let $\theta(s,t)=(-s^3-t^3,s^2t,st^2)$. Then the tangent cone
$\Z(f_3)$ is parameterized by 
$\theta$ as a map from $\mathbb{P}^1$ to $\mathbb{P}^2$.
When we need to compute the intersection
numbers between the rational curve $\Z(f_3)$ and the curve $\Z(f_4)$,
we can do that by studying the roots of the polynomial
$f_4(\theta)$. Expanding the polynomial gives
\begin{eqnarray*}
f_4(\theta)(s,t)&=&a_1s^{12}
-a_2s^{11}t+
(-a_3+a_4)s^{10}t^2
+(4a_1+a_5-a_7)s^9t^3\nonumber\\
&+&(-3a_2+a_6-a_8+a_{11})s^8t^4
+(-3a_3+2a_4-a_9+a_{12})s^7t^5\nonumber\\
&+&(6a_1+2a_5-a_7-a_{10}+a_{13})s^6t^6\\
&+&(-3a_2+2a_6-a_8+a_{14})s^5t^7
+(-3a_3+a_4-a_9+a_{15})s^4t^8\nonumber\\
&+&(4a_1+a_5-a_{10})s^3t^9
+(-a_2+a_6)s^2t^{10}
-a_3st^{11}
+a_1t^{12}\nonumber.
\end{eqnarray*}

This polynomial will have roots at $(0:1)$ and $(1:0)$ if and only
if $f_4(1,0,0)=a_1=0$.
When $a_1=0$ we may (by symmetry) assume $a_3\neq0$, so that $(0:1)$
is a simple root and $(1:0)$ is a root of multiplicity $m-1$.
Other roots of $f_4(\theta)$ correspond to intersections of $\Z(f_3)$
and $\Z(f_4)$ away from $(1:0:0)$. The multiplicity $m_i$ of each root
is equal to the corresponding intersection multiplicity, giving rise
to an $A_{m_i-1}$ singularity if $m_i>0$, as described by Proposition
\ref{prop_type_sing}, or a line $L_a\subset\Z(F)$ with $O$
as the only singular point if $m_i=1$.

The polynomial $f_4(\theta)$ defines a linear map from the
coefficient space $\kfield^{15}$ of $f_4$ to the space of homogeneous
polynomials of degree $12$ in $s$ and $t$. By elementary linear algebra,
we see that the image of this map is the set of polynomials of the
form
\begin{equation*}
b_0s^{12}+b_1s^{11}t+b_2s^{10}t^2+\dots+b_{12}t^{12}
\end{equation*}
where $b_0=b_{12}$. The kernel of the map corresponds to the set of
polynomials of the form $\ell f_3$ where $\ell$ is a linear form. This
means that $f_4(\theta)\equiv0$ if and only if $f_3$ is a
factor in $f_4$, making $\Z(F)$ reducible and not a monoid.

For every $m=0,2,3,4,\dots,12$ we can select $r$ parameter points
\begin{equation*}
p_1,\dots,p_r\in\mathbb{P}^1\setminus\{(0:1),(1:0)\}
\end{equation*}
 and positive
multiplicities $m_1,\dots,m_r$ with $m_1+\dots+m_r=12-m$ and try to
describe the polynomials $f_4$ such that $f_4(\theta)$ has a root of
multiplicity $m_i$ at $p_i$ for each $i=1,\dots,r$.

Still assuming $a_3\neq0$ whenever $a_1=0$, any such choice of
parameter points $p_1,\dots,p_r$ and multiplicities $m_1,\dots,m_r$
corresponds to a polynomial
$q=b_0s^{12}+b_1s^{11}t+\dots+b_{12}t^{12}$ that is, up to a nonzero
constant, uniquely determined.

Now, $q$ is equal to $f_4(\theta)$ for some $f_4$ if and only if
$b_0=b_{12}$. If $m\ge2$, then $q$
contains a factor $st^{m-1}$, so $b_0=b_{12}=0$,
giving $q=f_4(\theta)$ for some $f_4$.
In fact, when $m\ge2$ any choice of $p_1,\dots,p_r$ and $m_1,\dots,m_r$ 
with $m_1+\dots+m_r=12-m$ corresponds to a
four dimensional space of equations $f_4$ that gives this set of roots
and multiplicities in $f_4(\theta)$.  If $f_4^\prime$ is one such
$f_4$, then any other is of the form $\lambda f_4^\prime+\ell f_3$ for
some constant $\lambda\neq0$ and linear form $\ell$. All of these give
monoids that are projectively equivalent.

When $m=0$, we write $p_i=(\alpha_i:\beta_i)$ for $i=1,\dots,r$.
The condition $b_0=b_{12}$ on the coefficients of $q$ translates to
\begin{equation}\label{condition_points_NC}
\alpha_1^{m_1}\cdots\alpha_r^{m_r}=\beta_1^{m_1}\cdots\beta_r^{m_r}.
\end{equation}
This means that any choice of parameter points
$(\alpha_1:\beta_1),\dots,(\alpha_r:\beta_r)$ and multiplicities
$m_1,\dots,m_r$ with $m_1+\dots+m_r=12$ that satisfy
condition (\ref{condition_points_NC}) corresponds to a four
dimensional family $\lambda f_4^\prime+\ell f_3$, giving a unique
monoid up to projective equivalence.

For example, we can have an $A_{11}$ singularity only if $f_4(\theta)$
is of the form $(\alpha s-\beta t)^{12}$. Condition
(\ref{condition_points_NC}) implies that this can only happen for $12$
parameter points, all of the form $(1:\omega)$, where $\omega^{12}=1$.
Each such parameter point $(1:\omega)$
corresponds to a monoid
uniquely determined up to projective equivalence. However, since
there are six projective transformations of the plane that maps
$\Z(f_3)$ onto itself, this correspondence is not one to one. If
$\omega_1^{12}=\omega_2^{12}=1$, then $\omega_1$ and $\omega_2$ will
correspond to projectively equivalent monoids if and only if
$\omega_1^3=\omega_2^3$ or $\omega_1^3\omega_2^3=1$. This means that
there are three different quartic monoids with one $T_{3,3,4}$
singularity and one $A_{11}$ singularity. One corresponds to those
$\omega$ where $\omega^3=1$, one to those $\omega$ where $\omega^3=-1$,
and one to those $\omega$ where $\omega^6=-1$. The first two of these
have real representatives, $\omega=\pm1$.

It easy to see that for any set of multiplicities $m_1+\dots+m_r=12$,
we can find real points $p_1,\dots,p_r$ such that condition
(\ref{condition_points_NC}) is satisfied. This completely classifies
the possible configurations of singularities when $f_3$ is an
irreducible nodal curve.
\smallskip

{\bf Case 2.} The tangent cone is a cuspidal curve, and we can
assume $f_3(x_1,x_2,x_3)=x_1^3-x_2^2x_3$. The cuspidal curve is
singular at $(0:0:1)$ and can be parameterized by $\theta$ as a map
from $\mathbb{P}^1$ to $\mathbb{P}^2$ where
$\theta(s,t)=(s^2t,s^3,t^3)$. The intersection numbers are determined
by the degree $12$ polynomial $f_4(\theta)$. 
As in the previous case, $f_4(\theta)\equiv0$ if and
only if $f_3$ is a factor of $f_4$, and we will assume this is not
the case.
The multiplicity $m$ of the factor $s$ in $f_4(\theta)$ determines the
type of singularity at $O$. If $m=0$ (no factor $s$), then
$O$ is a $Q_{10}$ singularity. If $m=2$ or $m=3$, then $O$ is of type
$Q_{9+m}$.  If $m>3$, then $(0:0:1)$ is a singular point on $\Z(f_4)$,
so the monoid has a singular line and is not considered in this
article. Also, $m=1$ is not possible, since 
$f_{4}(\theta(s,t))=f_{4}(s^2t,s^3,t^3)$ cannot contain $st^{11}$ as a factor.

For each $m=0,2,3$ we can analyze the 
possible configurations %
of other singularities on the monoid. Similarly to the
previous case, any choice of parameter points
$p_1,\dots,p_r\in\mathbb{P}^1\setminus\{(0:1)\}$ and positive
multiplicities $m_1,\dots,m_r$ with $\sum m_i=12-m$ corresponds, up to
a nonzero constant, to a unique degree $12$ polynomial $q$.

When $m=2$ or $m=3$, for any choice of parameter values and associated
multiplicities, we can find a four dimensional family $f_4=\lambda
f_4^\prime+\ell f_3$ with the prescribed roots in $f_4(\theta)$. As
before, the family gives projectively equivalent monoids. 

When $m=0$,
one condition must be satisfied for $q$ to be of the form
$f_4(\theta)$, namely $b_{11}=0$, where $b_{11}$ is the coefficient
of $st^{11}$ in $q$.

For example, we can have an $A_{11}$ singularity only if $q$ is of the
form $(\alpha s-\beta t)^{12}$. The condition $b_{11}=0$ implies that
either $q=\lambda s^{12}$ or $q=\lambda t^{12}$. The first case gives
a surface with a singular line, while the other gives a monoid with an
$A_{11}$ singularity (see Figure \ref{fig_ex_2}). The line from $O$ to the $A_{11}$ singularity
corresponds to the inflection point of $\Z(f_3)$.

For any set of multiplicities $m_1,\dots,m_r$ with
$m_1+\dots+m_r=12$, it is not hard to see that there exist real points
$p_1,\dots,p_r$ such that the condition $b_{11}=0$ is satisfied. It suffices to
take $p_{i}=(\alpha_{i}:1)$, with $\sum m_{i}\alpha_{i}=0$ (the condition 
corresponding to $b_{11}=0$).
This completely classifies the possible configurations of singularities
when $f_3$ is a cuspidal curve.
\smallskip

{\bf Case 3.} The tangent cone is the product of a conic and a
line that is not tangent to the conic, and
we can assume $f_3=x_3(x_1x_2+x_3^2)$. Then $\Z(f_3)$ is singular at
$(1:0:0)$ and $(0:1:0)$, the intersections of the conic
$\Z(x_1x_2+x_3^2)$ and the line $\Z(x_3)$. For each $f_4$ we can
associate four integers:
\begin{equation*}
\begin{tabular}{ll}
$j_0:=\I_{(1:0:0)}(x_1x_2+x_3^2,f_4)$, & $k_0:=\I_{(1:0:0)}(x_3,f_4)$, \\
$j_1:=\I_{(0:1:0)}(x_1x_2+x_3^2,f_4)$, & $k_1:=\I_{(0:1:0)}(x_3,f_4)$.
\end{tabular}
\end{equation*}

We see that $k_0>0\Leftrightarrow f_4(1:0:0)=0\Leftrightarrow j_0>0$,
and that $\Z(f_4)$ is singular at $(1:0:0)$ if and only if $k_0$ and
$j_0$ both are bigger than one.  These cases imply a singular line on
the monoid, and are not considered in this article.  The same holds for
$k_1$, $j_1$ and the point $(0:1:0)$.

Define $r_i=\max(j_i,k_i)$ for $i=1,2$. Then, by \cite{MR693454}, $O$ will be a
singularity of type $T_{3,4+r_0,4+r_1}$ if $r_0\le r_1$, or of type 
$T_{3,4+r_1,4+r_0}$ if $r_0\ge r_1$.

We can parameterize the line $\Z(x_3)$ by 
$\theta_1$ where $\theta_1(s,t)=(s,t,0)$, and the conic
$\Z(x_1x_2+x_3^2)$ by $\theta_2$ where $\theta_2(s,t)=(s^2,-t^2,st)$.
Similarly to the previous cases, roots of $f_4(\theta_1)$ correspond to
intersections between $\Z(f_4)$ and the line $\Z(x_3)$, while roots of
$f_4(\theta_2)$ correspond to intersections between $\Z(f_4)$ and the
conic $\Z(x_1x_3+x_3^2)$.

For any legal values of of $j_0$, $j_1$, $k_0$ and $k_1$, parameter points
\begin{equation*}
(\alpha_1:\beta_1),\dots,(\alpha_{m_r}:\beta_{m_r})\in
\mathbb{P}^1\setminus\{(0:1),(1:0)\},
\end{equation*}
with multiplicities $m_1,\dots,m_r$ such that
$m_1+\dots+m_r=4-k_0-k_1$, and parameter points
\begin{equation*}
(\alpha_1^\prime:\beta_1^\prime),\dots,(\alpha_{m_r^\prime}^\prime:
\beta_{m_r^\prime}^\prime)\in\mathbb{P}^1\setminus\{(0:1),(1:0)\},
\end{equation*}
with multiplicities $m_1^\prime,\dots,m_{r^\prime}^\prime$ such that
$m_1^\prime+\dots+m_{r^\prime}^\prime=8-j_0-j_1$, we can fix
polynomials $q_1$ and $q_2$ such that
\begin{itemize}
\item[$\bullet$] $q_1$ is nonzero, of degree $4$, and has factors $s^{k_1}$, $t^{k_0}$ and
  $(\beta_is-\alpha_it)^{m_i}$ for $i=1,\dots,r$,
\item[$\bullet$] $q_2$ is nonzero, of degree $8$, and has factors $s^{j_1}$, $t^{j_0}$ and 
$(\beta_i^\prime s-\alpha_i^\prime t)^{m_i^\prime}$ for $i=1,\dots,r^\prime$.
\end{itemize}
Now $q_1$ and $q_2$ are determined up to multiplication by nonzero
constants. Write
$q_1=b_0s^4+\dots+b_4t^4$ and $q_2=c_0s^8+\dots+c_8t^8$.

The classification of singularities on the monoid consists of
describing the conditions on the parameter points and nonzero
constants $\lambda_1$ and $\lambda_2$ for the pair
$(\lambda_1q_1,\lambda_2q_2)$ to be on the form
$(f_4(\theta_1),f_4(\theta_2))$ for some $f_4$.

Similarly to the previous cases, $f_4(\theta_1)\equiv0$ if and only if
$x_3$ is a factor in $f_4$ and $f_4(\theta_2)\equiv0$ if and only if
$x_1x_2+x_3^2$ is a factor in $f_4$. Since $f_3=x_3(x_1x_2+x_3^2)$,
both cases will make the monoid reducible, so
we only consider $\lambda_1,\lambda_2\neq0$.

We use linear algebra to study the relationship between the
coefficients $a_1\dots a_{15}$ of $f_4$ and the polynomials $q_1$ and
$q_2$. We find $(\lambda_1q_1,\lambda_2q_2)$ to be of the form
$(f_4(\theta_1),f_4(\theta_2))$ if and only if
$\lambda_1b_0=\lambda_2c_0$ and $\lambda_1b_4=\lambda_2c_8$.
Furthermore, the pair $(\lambda_1q_1,\lambda_2q_2)$ will fix $f_4$
modulo $f_3$.
Since $f_4$ and $\lambda f_4$ correspond to projectively equivalent
monoids for any $\lambda\neq0$, it is the ratio $\lambda_1/\lambda_2$,
and not $\lambda_1$ and $\lambda_2$, that is important.

Recall that $k_0>0\Leftrightarrow j_0>0$ and $k_1>0\Leftrightarrow
j_1>0$.  If $k_0>0$ and $k_1>0$, then $b_0=c_0=b_4=c_8=0$, so for any
$\lambda_1,\lambda_2\neq0$ we have
$(\lambda_1q_1,\lambda_2q_2)=(f_4(\theta_1),f_4(\theta_2))$ for some
$f_4$.
Varying $\lambda_1/\lambda_2$ will give a one-parameter family of
monoids for each choice of multiplicities and parameter
points.

If $k_0=0$ and $k_1>0$, then $b_0=c_0=0$. The condition
$\lambda_1b_4=\lambda_2c_8$ implies
$\lambda_1/\lambda_2=c_8/b_4$.
This means that any choice of multiplicities and parameter points will
give a unique monoid up to projective equivalence. The same goes for
the case where $k_0>0$ and $k_1=0$.

Finally, consider the case where $k_0=k_1=0$. For
$(\lambda_1q_1,\lambda_2q_2)$ to be of the form
$(f_4(\theta_1),f_4(\theta_2))$ we must have
$\lambda_1/\lambda_2=c_8/b_4=c_0/b_0$. This translates into a
condition on the parameter points, namely
\begin{equation}\label{condition_points_QC}
\frac{(\beta_1^\prime)^{m_1^\prime}\cdots(\beta_{r\prime}^\prime)^{m_{r\prime}^\prime}}
{\beta_1^{m_1}\cdots\beta_r^{m_r}}=
\frac{(\alpha_1^\prime)^{m_1^\prime}\cdots(\alpha_{r^\prime}^\prime)^{m_{r^\prime}^\prime}}
{\alpha_1^{m_1}\cdots\alpha_r^{m_r}}.
\end{equation}
In other words, if condition (\ref{condition_points_QC}) holds, we
have a unique monoid up to projective equivalence.

It is easy to see that for any choice of multiplicities, it is
possible to find real parameter points such that condition
(\ref{condition_points_QC}) is satisfied. This completes the
classification of possible singularities when the tangent cone is a conic plus
a chordal line.
\smallskip

{\bf Case 4.} The tangent cone is the product of a conic and a
line tangent to the conic, and we can assume $f_3=x_3(x_1x_3+x_2^2)$.
Now $\Z(f_3)$ is singular at $(1:0:0)$. For each $f_4$ we can
associate two integers
\begin{equation*}
j_0:=\I_{(1:0:0)}(x_1x_3+x_2^2,f_4)\hspace{2em}\text{and}\hspace{2em}k_0:=\I_{(1:0:0)}(x_3,f_4).
\end{equation*}
We have $j_0>0\Leftrightarrow k_0>0$, $j_0>1\Leftrightarrow k_0>1$.
Furthermore, $j_0$ and $k_0$ are both greater than $2$ if and only if
$\Z(f_4)$ is singular at $(1:0:0)$, a case we have excluded.
The singularity at $O$ will be of the $S$ series, from \cite{ArnoldLondon1}, 
\cite{ArnoldLondon2}.

We can parameterize the conic $\Z(x_1x_3+x_2^2)$ by $\theta_2$ and the
line $\Z(x_3)$ by $\theta_1$ where $\theta_2(s:t)=(s^2:st:-t^2)$ and
$\theta_1(s:t)=(s:t:0)$.
As in the previous case, the monoid is reducible if
and only if $f_4(\theta_1)\equiv0$ or $f_4(\theta_2)\equiv0$.
Consider two nonzero polynomials
\begin{eqnarray*}
q_1&=&b_0s^4+b_1s^3t+b_2s^2t^2+b_3st^3+b_4t^4 \\
q_2&=&c_0s^8+c_1s^7t+\dots+c_7st^7+c_8t^8.
\end{eqnarray*}
Now $(\lambda_1q_1,\lambda_2q_2)=(f_4(\theta_1),f_4(\theta_2))$ for
some $f_4$ if and only if $\lambda_1b_0=\lambda_2c_0$ and
$\lambda_1b_1=\lambda_2c_1$. As before, only the
  cases where $\lambda_1,\lambda_2\neq0$ are interesting. We see that
  $(\lambda_1q_1,\lambda_2q_2)=(f_4(\theta_1),f_4(\theta_2))$ for some
  $\lambda_1,\lambda_2\neq0$ if and only if the following hold:
\begin{itemize}
\item[$\bullet$] $b_0=0\leftrightarrow c_0=0$
and $b_1=0\leftrightarrow c_1=0$
\item[$\bullet$] $b_0c_1=b_1c_0$.
\end{itemize}

The classification of other singularities (than $O$) is very similar to
the previous case. Roots of $f_4(\theta_1)$ and $f_4(\theta_2)$ away
from $(1:0)$ correspond to intersections of $\Z(f_3)$ and $\Z(f_4)$
away from the singular point of $\Z(f_3)$, and when one such
intersection is multiple, there is a corresponding 
singularity on the monoid.

Now assume $(\lambda_1q_1,\lambda_2q_2)=(f_4(\theta_1),f_4(\theta_2))$
for some $\lambda_1,\lambda_2\neq0$ and some $f_4$.
If $b_0\neq0$
(equivalent to $c_0\neq0$) then $j_0=k_0=0$ and
$\lambda_1/\lambda_2=c_0/b_0$.  If $b_0=c_0=0$ and $b_1\neq0$
(equivalent to $c_1\neq0$), then $j_0=k_0=1$, and
$\lambda_1/\lambda_2=c_1/b_1$. If $b_0=b_1=c_0=c_1=0$, then
$j_0,k_0>1$ and any value of $\lambda_1/\lambda_2$ will give
$(\lambda_1q_1,\lambda_2q_2)$ of the form
$(f_4(\theta_1),f_4(\theta_2))$ for some $f_4$. Thus we get a 
one-dimensional family of monoids for this choice of $q_1$ and $q_2$.

Now consider the possible configurations of other singularities on the
monoid.  Assume that $j_0^\prime\le8$ and $k_0^\prime\le4$ are
nonnegative integers such that $j_0>0\leftrightarrow k_0>0$ and
$j_0>1\leftrightarrow k_0>1$.
For any set of multiplicities
$m_1,\dots,m_r$ with $m_1+\dots+m_r=4-k_0^\prime$ and
$m_1^\prime,\dots,m_{r^\prime}^\prime$ with
$m_1^\prime+\dots+m_{r^\prime}^\prime=8-j_0^\prime$, there exists a
polynomial $f_4$ with real coefficients such that $f_4(\theta_1)$ has
real roots away from $(1:0)$ with multiplicities $m_1,\dots,m_r$, and 
$f_4(\theta_2)$ has real roots away from $(1:0)$ with multiplicities
$m_1^\prime,\dots,m_{r^\prime}^\prime$. Furthermore, for this $f_4$ we have
$k_0=k_0^\prime$ and $j_0=j_0^\prime$. Proposition \ref{prop_type_sing}
will give the singularities that occur in addition to $O$.

This completes the classification of the singularities on a quartic monoid (other than $O$)
when the tangent cone is a conic plus a tangent.
\smallskip

{\bf Case 5.} The tangent cone is three general lines, and we
assume $f_3=x_1x_2x_3$.

For each $f_4$ we associate six integers,
\begin{eqnarray*}
k_2:=\I_{(1:0:0)}(f_4,x_2),
& l_1:=\I_{(0:1:0)}(f_4,x_1),
& m_1:=\I_{(0:0:1)}(f_4,x_1),\\
k_3:=\I_{(1:0:0)}(f_4,x_3),
& l_3:=\I_{(0:1:0)}(f_4,x_3),
& m_2:=\I_{(0:0:1)}(f_4,x_2).
\end{eqnarray*}
Now $k_2>0\Leftrightarrow k_3>0$, $l_1>0\Leftrightarrow l_3>0$, and
$m_1>0\Leftrightarrow m_2>0$. If both $k_2$ and $k_3$ are greater than
$1$, then the monoid has a singular line, a case we have excluded. 
The same goes for the pairs $(l_1,l_3)$ and $(m_1,m_2)$.

When the monoid does not have a singular line, we define
$j_k=\max(k_2,k_3)$, $j_l=\max(l_1,l_3)$ and $j_m=\max(m_1,m_2)$. If
$j_k\le j_l\le j_m$, then \cite{MR693454} gives that $O$ is a
$T_{4+j_k,4+j_l,4+j_m}$ singularity.

The three lines $\Z(x_1)$, $\Z(x_2)$ and $\Z(x_3)$ are parameterized
by $\theta_1$, $\theta_2$ and $\theta_3$ where
$\theta_1(s:t)=(0:s:t)$, $\theta_2(s:t)=(s:0:t)$ and
$\theta_3(s:t)=(s:t:0)$. Roots of the polynomial $f_4(\theta_i)$
away from $(1:0)$ and $(0:1)$ correspond to intersections between
$\Z(f_4)$ and $\Z(x_i)$ away from the singular points of $\Z(f_3)$.

As before, we are only interested in the cases where
none of $f_4(\theta_i)\equiv0$ for $i=1,2,3$, as this would make the monoid
reducible.

For the study of other singularities on the monoid we consider nonzero
polynomials
\begin{eqnarray*}
q_1&=&b_0s^4+b_1s^3t+b_2s^2t^2+b_3st^3+b_4t^4,\\
q_2&=&c_0s^4+c_1s^3t+c_2s^2t^2+c_3st^3+c_4t^4,\\
q_3&=&d_0s^4+d_1s^3t+d_2s^2t^2+d_3st^3+d_4t^4.
\end{eqnarray*}
Linear algebra shows that
$(\lambda_1q_1,\lambda_2q_2,\lambda_3q_3)=(f_4(\theta_1),f_4(\theta_2),f_4(\theta_3))$
for some $f_4$ if and only if $\lambda_1b_0=\lambda_3d_4$,
$\lambda_1b_4=\lambda_2c_4$, and $\lambda_2c_0=\lambda_3d_0$.
A simple analysis shows the following: There exist
$\lambda_1,\lambda_2,\lambda_3\neq0$ such that
\[
(\lambda_1q_1,\lambda_2q_2,\lambda_3q_3)
=(f_4(\theta_1),f_4(\theta_2),f_4(\theta_3))\]
for some $f_4$, and such that $\Z(f_4)$ and $\Z(f_3)$ have no common
singular point if and only if all of the following hold:
\begin{itemize}
\item[$\bullet$] $b_0=0\leftrightarrow d_4=0$ and $b_0=d_4=0\rightarrow (b_1\neq0$ or $d_3\neq0)$,
\item[$\bullet$] $b_4=0\leftrightarrow c_4=0$ and $b_4=c_4=0\rightarrow (b_3\neq0$ or $c_3\neq0)$,
\item[$\bullet$] $c_0=4\leftrightarrow d_0=0$ and $c_0=d_0=0\rightarrow (c_1\neq0$ or $d_1\neq0)$,
\item[$\bullet$] $b_0c_4d_0=b_4c_0d_4$.
\end{itemize}

Similarly to the previous cases we can classify the possible
configurations of other singularities by varying the multiplicities of
the roots of the polynomials $q_1$, $q_2$ and $q_3$. Only the
multiplicities of the roots $(0:1)$ and $(1:0)$ affect the first three
bullet points above. Then, for any set of multiplicities of the
rest of the roots, we can find $q_1$, $q_2$ and $q_3$ such that 
the last bullet point is satisfied. This completes
the classification when $\Z(f_3)$ is the product of three general lines.
\smallskip

{\bf Case 6.} The tangent cone is three lines meeting in a point,
and we can assume that $f_3=x_2^3-x_2x_3^2$. We write
$f_3=\ell_1\ell_2\ell_3$ where $\ell_1=x_2$, $\ell_2=x_2-x_3$ and
$\ell_3=x_2+x_3$, representing the three lines going through the
singular point $(1:0:0)$. For each $f_4$ we associate three integers
$j_1$, $j_2$ and $j_3$ defined as the intersection numbers
$j_i=\I_{(1:0:0)}(f_4,\ell_i)$. We see that $j_1=0\Leftrightarrow
j_2=0\Leftrightarrow j_3=0$, and that $\Z(f_4)$ is singular at
$(1:0:0)$ if and only if two of the integers $j_1$, $j_2$, $j_3$ are
greater then one. (Then all of them will be greater than one.)
The singularity will be of the $U$ series \cite{ArnoldLondon1}, \cite{ArnoldLondon2}.

The three lines $\Z(\ell_1)$, $\Z(\ell_2)$ and $\Z(\ell_3)$ can be
parameterized by $\theta_1$, $\theta_2$, and $\theta_3$ where 
$\theta_1(s:t)=(s:0:t)$, $\theta_2(s:t)=(s:t:t)$ and $\theta_2(s:t)=(s:t:-t)$.

For the study of other singularities on the monoid we consider nonzero
polynomials
\begin{eqnarray*}
q_1&=&b_0s^4+b_1s^3t+b_2s^2t^2+b_3st^3+b_4t^4,\\
q_2&=&c_0s^4+c_1s^3t+c_2s^2t^2+c_3st^3+c_4t^4,\\
q_3&=&d_0s^4+d_1s^3t+d_2s^2t^2+d_3st^3+d_4t^4.
\end{eqnarray*}
Linear algebra shows that
$(\lambda_1q_1,\lambda_2q_2,\lambda_3q_3)=(f_4(\theta_1),f_4(\theta_2),f_4(\theta_3))$
for some $f_4$ if and only if
$\lambda_1b_0=\lambda_2c_4=\lambda_3d_0$, and
$2\lambda_1b_1=\lambda_2c_1+\lambda_3d_1$. %
There exist $\lambda_1,\lambda_2,\lambda_3\neq0$ such that 
$(\lambda_1q_1,\lambda_2q_2,\lambda_3q_3)=(f_4(\theta_1),f_4(\theta_2),f_4(\theta_3))$
for some $f_4$ and such that $\Z(f_4)$ and $\Z(f_3)$  have no common
singular point if and only if all of the following hold:
\begin{itemize}
\item[$\bullet$] $b_0=0\leftrightarrow c_0=0\leftrightarrow d_0=0$,
\item[$\bullet$]
  if $b_0=c_0=d_0=0$, then at least two of $b_1$, $c_1$, and $d_1$
  are different from zero,
\item[$\bullet$] $2b_1c_0d_0=b_0c_1d_0+b_0c_0d_1$.
\end{itemize}

As in all the previous cases we can classify the possible
configurations of other singularities for all possible $j_1,j_2,j_3$.
As before, the first bullet point only affect the
multiplicity of the factor $t$ in $q_1$, $q_2$ and $q_3$. For any set
of multiplicities for the rest of the roots, we can find $q_1$, $q_2$,
$q_3$ with real roots of the given multiplicities such that the last bullet point is
satisfied. This completes the classification of the singularities (other than $O$)
when $\Z(f_3)$ is three lines meeting in a point.
\smallskip

{\bf Case 7.} The tangent cone is a double line plus a line, and
we can assume $f_3=x_2x_3^2$. The tangent cone is singular along the
line $\Z(x_3)$. The line $\Z(x_2)$ is parameterized by $\theta_1$ and
the line $\Z(x_3)$ is parameterized by $\theta_2$ where
$\theta_1(s:t)=(s:0:t)$ and $\theta_2(s:t)=(s:t:0)$.
The monoid is reducible if and only if $f_4(\theta_1)$ or
$f_4(\theta_2)$ is identically zero, so we assume that neither is
identically zero.
For each $f_4$ we associate two integers, $j_0:=\I_{(1:0:0)}(f_4,x_2)$
and $k_0:=\I_{(1:0:0)}(f_4,x_3)$. Furthermore, we write $f_4(\theta_2)$ as
a product of linear factors
\begin{equation*}
f_4(\theta_2)=\lambda s^{k_0}\prod_{i=0}^r(\alpha_is-t)^{m_i}.
\end{equation*}
Now the singularity at $O$ will be of the $V$ series and depends
on $j_0$, $k_0$ and $m_1,\dots,m_r$.

Other singularities on the monoid correspond to intersections of
$\Z(f_4)$ and the line $\Z(x_2)$ away from $(1:0:0)$. Each such
intersection corresponds to a root in the polynomial $f_4(\theta_1)$
different from $(1:0)$. Let $j_0^\prime\le4$ and $k_0^\prime\le4$
be integers such that $j_0>0\leftrightarrow k_0>0$. Then, for any
homogeneous polynomials $q_1$, $q_2$ in $s,t$ of degree 4 such that
$s$ is a factor of multiplicity $j_0^\prime$ in $q_1$ and of multiplicity
$k_0^\prime$ in $q_2$, there is a polynomial $f_4$ and nonzero constants
$\lambda_1$ and $\lambda_2$ such that $k_0=k_0^\prime$,
$j_0=j_0^\prime$ and
$(\lambda_1q_1,\lambda_2q_2)=(f_4(\theta_1),f_4(\theta_2))$.
Furthermore, if $q_1$ and $q_2$ have real coefficients, then $f_4$ can
be selected with real coefficients. This follows from an analysis
similar to case 5 and completes the classification of singularities
when the tangent cone is a product of a line and a double line.

{\bf Case 8.} The tangent cone is a triple line, and we assume
that $f_3=x_3^3$. The line $\Z(x_3)$ is parameterized by $\theta$
where $\theta(s,t)=(s,t,0)$. Assume that the polynomial $f_4(\theta)$
has $r$ distinct roots with multiplicities $m_1,\dots,m_r$. (As before
$f_4(\theta)\equiv 0$ if and only if the monoid is reducible.) Then the
type of the singularity at $O$ will be of the $V^\prime$ series \cite[p.~267]{MR777682}.
The integers $m_1,\dots,m_r$ are constant under right equivalence
over $\CC$. Note that one can construct examples of monoids that
are right equivalent over $\CC$, but not over $\RR$ (see Figure \ref{fig_V}).

\begin{figure}[!htbp]
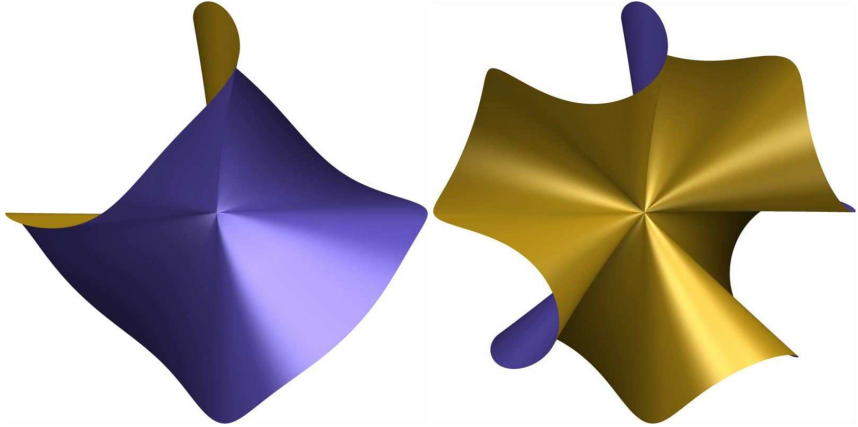
%
\includegraphics[width=.46\textwidth]{fig_z3+xy3+x3y.epsf}
\includegraphics[width=.46\textwidth]{fig_z3+xy3-x3y.epsf}%
\caption{The monoids $\Z(z^3+xy^3+x^3y)$ and
$\Z(z^3+xy^3-x^3y)$ are right equivalent over $\CC$ but not over
$\RR$.
}%
\label{fig_V}%
\end{figure}

The tangent cone is singular everywhere, so there can be no
other singularities on the monoid.
\smallskip

{\bf Case 9.} The tangent cone is a smooth cubic curve, and we write
$f_3=x_1^3+x_2^3+x_3^3+3ax_1x_2x_3$ where $a^3\neq-1$. This is
a one-parameter family of elliptic curves, so we cannot use the
parameterization technique of the other cases.  The singularity at $O$ 
will be a $P_8$ singularity (cf. \cite[p.~185]{MR777682}), and other singularities
correspond to intersections between $\Z(f_3)$ and $\Z(f_4)$, as described by
Proposition \ref{prop_type_sing}.

To classify the possible configurations of singularities on a monoid
with a nonsingular (projective) tangent cone, we need to answer the following
question: For any positive integers $m_1,\dots,m_r$ such that
$\sum_{i=1}^rm_i=12$, does there, for \emph{some}
$a\in{\RR}\setminus\{-1\}$, exist a polynomial $f_4$ with real
coefficients such that $\Z(f_3,f_4)=\{p_1,\dots,p_r\}\in\mathbb{P}^2({\RR})$ and
$\I_{p_i}(f_3,f_4)=m_i$ for $i=1,\dots,r$? Rohn \cite[p.~63]{MR1510280}
says that one can always find curves $\Z(f_3)$, $\Z(f_4)$ with this property. Here we 
shall show that
for \emph{any}  $a\in{\RR}\setminus\{-1\}$ we can find a suitable $f_{4}$.

In fact, in almost all cases $f_4$ can be constructed as a product of linear
and quadratic terms in a simple way.
The difficult cases are $(m_1,m_2)=(11,1)$, $(m_1,m_2,m_3)=(8,3,1)$, and
$(m_1,m_2)=(5,7)$. 
For example, the case where $(m_1,m_2,m_3)=(3,4,5)$
can be constructed as follows: Let $f_4=\ell_1\ell_2\ell_3^2$ where
$\ell_1$ and $\ell_2$ define tangent lines at inflection points $p_1$
and $p_3$ of $\Z(f_3)$.
Let $\ell_3$ define a line that intersects $\Z(f_3)$ once at $p_3$ and
twice at another point $p_2$. Note that the points $p_1$, $p_2$ and
$p_3$ can be found for any $a\in\RR\setminus\{-1\}$.

The case $(m_1,m_2)=(11,1)$ is also possible for every
$a\in\RR\setminus\{-1\}$. For any point $p$ on $\Z(f_3)$ there
exists an $f_4$ such that $\I_p(f_3,f_4)\ge11$. For all except a
finite number of points, we have equality \cite{MR2019943}, so the case
$(m_1,m_2)=(11,1)$ is possible for any
$a\in\RR\setminus\{-1\}$. The case $(m_1,m_2,m_3)=(8,3,1)$ is
similar, but we need to let $f_4$ be a product of the tangent at an
inflection point with another cubic.

The case $(m_1,m_2)=(5,7)$ is harder. Let $a=0$. Then we can construct
a conic $C$ that intersects $\Z(f_3)$ with multiplicity five in one point
and multiplicity one in an inflection point, and choosing $\Z(f_4)$ as the union
of $C$ and twice the tangent line through the inflection point will give the
desired example. The same can be done for $a=-4/3$. By using the computer
algebra system {\sc Singular} \cite{GPS05}
we can show that these constructions can be continuously extended to any
$a\in\RR\setminus\{-1\}$.  This completes the classification of
singularities on a monoid when the tangent cone is smooth.
\medskip

In the Cases $3$, $5$, and $6$, not all real equations of a given type can
be transformed to the chosen forms by a real transformation.

In Case
$3$ the conic may not intersect the line in two real points, but rather
in two complex conjugate points. Then we can assume
$f_3=x_3(x_1x_3+x_1^2+x_2^2)$, and the singular points are $(1:\pm
i:0)$. For any real $f_4$, we must have
\[
\I_{(1:i:0)}(x_1x_3+x_1^2+x_2^2,f_4)=\I_{(1:-i:0)}(x_1x_3+x_1^2+x_2^2,f_4)\]
and 
\[
\I_{(1:i:0)}(x_3,f_4)=\I_{(1:-i:0)}(x_3,f_4),\] 
so only the cases
where $j_0=j_1$ and $k_0=k_1$ are possible. Apart from that, no other
restrictions apply.

In Case $5$, two of the lines can be complex conjugate, and we assume
$f_3=x_3(x_1^2+x_2^2)$. A configuration from the previous analysis is
possible for real coefficients of $f_4$ if and only if $m_1=m_2$,
$k_2=l_1$, and $k_3=l_3$. Furthermore, only the singularities 
that correspond to the line $\Z(x_3)$ will be real.

In Case $6$, two of the lines can be complex conjugate, and then we
may assume $f_3=x_2^3+x_3^3$. Now, if $j_{3}$ denotes the intersection number of 
$\Z(f_{4})$
with the real line $\Z(x_2+x_3)$,
precisely the cases where $j_1=j_2$
are possible, and only intersections with the line $\Z(x_2+x_3)$ may
contribute to real singularities.

This concludes the classification of real and complex
singularities on real monoids of degree $4$.
\end{proof}

\textbf{Remark.} In order to describe the various monoid singularities, Rohn 
\cite{MR1510280} computes
the ``class reduction" due to the presence of the singularity, in (almost) all cases. 
(The class is the degree of the dual surface 
\cite[p.~262]{MR510551}.)
The class reduction is equal to the local intersection multiplicity of the surface with 
two general polar surfaces. This
intersection multiplicity is equal to the sum of the Milnor number and the Milnor number of
a general plane section through the singular point \cite[Cor. 1.5, p.~320]{MR0374482}. 
It is not hard to see that a general plane
section has either a $D_{4}$ (Cases 1--6, 9), $D_{5}$ (Case 7), or $E_{6}$ (Case 8) singularity. Therefore one can 
retrieve the Milnor number of each monoid singularity from Rohn's work.

\subsection*{Acknowledgements}
We would like to thank the referees for helpful comments.
This research was supported by the European Union through the project IST 2001--35512 
`Intersection algorithms for geometry based IT applications using approximate algebraic
methods` (GAIA II).


\bibliographystyle{plain}
\bibliography{biblio}

\end{document}